\documentclass[12pt]{article}
\usepackage{amsmath,amsfonts,amsthm,amssymb,eucal}
\begin{document}
\newtheorem{lem}{Lemma}
\newtheorem{teo}{Theorem}
\newtheorem{prop}{Proposition}
\pagestyle{plain}
\title{Non-Archimedean Shift Operators}
\author{Anatoly N. Kochubei\\
\footnotesize Institute of Mathematics,\\
\footnotesize National Academy of Sciences of Ukraine,\\
\footnotesize Tereshchenkivska 3, Kiev, 01601 Ukraine
\\ \footnotesize E-mail: \ kochubei@i.com.ua}
\date{}
\maketitle

\bigskip
\begin{abstract}
We introduce and study non-Archimedean analogs of the operators of unilateral shift and backward shift playing crucial roles in the classical theory of nonselfadjoint operators. In particular, we find various functional models of these operators having both common and different features compared to their classical counterparts.
\end{abstract}

\bigskip
{\bf Key words: }unilateral shift; backward shift; non-Archimedean Banach space; Tate algebra; indefinite sum; annihilation operator

{\bf MSC 2010}. Primary: 47S10. Secondary: 47B37.

\bigskip
\section{INTRODUCTION}
In the recent non-Archimedean spectral theorem \cite{K2}, the author found a class of linear operators on non-Archimedean Banach spaces, which admit spectral decompositions resembling the ones known for classical normal operators. Following the development of classical operator theory (see \cite{H,N}), in this note we introduce some model non-normal operators, the non-Archimedean counterparts of the unilateral shift and backward shift operators. We describe several of their functional models expressed in terms of basic notions and constructions of $p$-adic analysis. Some properties of these operators are parallel to the classical patterns, others are quite different. For example, the lattice of invariant subspaces of the non-Archimedean  unilateral shift is indexed by polynomials (in the classical Beurling theorem this role is played by inner functions), the operator itself can be seen also as an analog of the Volterra integration operator whose properties are (classically) very far from those of the unilateral shift. The author hopes that the examples given in this paper will become the first steps of a still non-existent theory of non-Archimedean non-normal operators.

Below we use the standard notation of $p$-adic analysis \cite{PS,R,S}. For a prime number $p$, $\mathbb Q_p$ is the field of $p$-adic numbers endowed with the absolute value $|\cdot |_p$, $\mathbb Z_P$ is the ring of integers (= the unit ball) in $\mathbb Q_p$, $\mathbb C_p$ is the completion, with respect to the extension of the absolute value $|\cdot |_p$, of an algebraic closure of $\mathbb Q_p$, $\mathbb A_p$ is the closed unit ball in $\mathbb C_p$.

We will use the following Banach spaces over $\mathbb C_p$. The space $c_0$ consists of all the sequences $(x_0,x_1,x_2,\ldots )$ of elements from $\mathbb C_p$, such that $x_n\to 0$ in $\mathbb C_p$, as $n\to \infty$. The space $l^\infty$, conjugate to $c_0$, consists of all bounded sequences. The norms on both spaces are equal to $\sup\limits_n|x_n|_p$. The space of all continuous functions on $\mathbb Z_p$ with values in $\mathbb C_p$, with the supremum norm, will be denoted by $C(\mathbb Z_p,\mathbb C_p)$. The space $H(\mathbb A_p)$ consists of all analytic elements on $\mathbb A_p$, that is the series $f(z)=\sum\limits_{n=0}^\infty a_nz^n$, $z\in \mathbb A_p$, where $(a_n)\in c_0$. The norm on $H(\mathbb A_p)$ is
$$
\|f\|=\sup\limits_{z\in \mathbb A_p}|f(z)|_p=\sup\limits_{n\ge 0}|a_n|_p
$$
(see Theorem 6.4.3 in \cite{R}). The space $H(\mathbb A_p)$ is actually a commutative Banach algebra (often called the Tate algebra).

If a non-Archimedean Banach space $X$ over $\mathbb C_p$ possesses an orthonormal basis (see \cite{PS,R,S} for the definition of the latter in the non-Archimedean case), then the expansion with respect to the basis establishes an isometric isomorphism of $X$ and $c_0$. Therefore the role of $c_0$ is similar to that of $l^2$ in classical functional analysis, while $C(\mathbb Z_p,\mathbb C_p)$ and $H(\mathbb A_p)$ can be seen as counterparts of the Hardy space $H^2$.

\section{Unilateral Shift}

The unilateral shift operator $S:\ c_0\to c_0$ is defined as
\begin{equation}
S(x_0,x_1,x_2,\ldots )=(0,x_0,x_1,x_2,\ldots ),
\end{equation}
for each sequence $(x_n)\in c_0$. We give some functional models of the operator $S$.

\medskip
1) Define an operator $S_1$ on $H(\mathbb A_p)$ setting $(S_1f)(z)=zf(z)$, $z\in \mathbb A_p$, for $f(z)=\sum\limits_{n=0}^\infty a_nz^n$, $(a_n)\in c_0$. We have
$$
(S_1f)(z)=\sum\limits_{n=1}^\infty a_{n-1}z^n=\sum\limits_{n=0}^\infty (Sa)_nz^n
$$
where $Sa=S(a_0,a_1,a_2,\ldots )=(0,a_0,a_1,a_2,\ldots )$, so that under the isomorphism $H(\mathbb A_p)\to c_0$, $f\mapsto (a_0,a_1,a_2,\ldots )$, the operator $S_1$ is transformed into $S$. Note that a formally similar multiplication operator on $C(\mathbb Z_p,\mathbb C_p)$ has a different sequence space representation; in fact, it is normal \cite{K2}, thus its properties are, as we will see below, quite different from those of $S_1$.

2) Let $S_2$ be an operator on $C(\mathbb Z_p,\mathbb C_p)$ defined as follows. For $\varphi \in C(\mathbb Z_p,\mathbb C_p)$, set
$$
(S_2\varphi )(n)=\sum\limits_{j=0}^{n-1}f(j),\quad n\in \mathbb N.
$$
The function $S_2\varphi$ (called the indefinite sum of $\varphi$) extends continuously onto $\mathbb Z_p$ \cite{R,S}.

It is well known that the Mahler polynomials
$$
P_n(x)=\frac{x(x-1)\cdots (x-n+1)}{n!},\quad n\ge 1;\quad P_0(x)=1,
$$
form an orthonormal basis in $C(\mathbb Z_p,\mathbb C_p)$. If $\varphi =\sum\limits_{n=0}^\infty b_nP_n$, $(b_n)\in c_0$, we have \cite{R,S}
$$
S_2\varphi =\sum\limits_{n=1}^\infty b_{n-1}P_n=\sum\limits_{n=0}^\infty (Sb)_nP_n.
$$

The operator $S_2$ can be written also in terms of the shifted convolution $\varphi \underline{*}\psi$ of functions $\varphi ,\psi \in C(\mathbb Z_p,\mathbb C_p)$. By definition,
$$
(\varphi \underline{*}\psi )(n)=\sum\limits_{i+j=n-1}\varphi (i)\psi (j),\quad n\in \mathbb N,
$$
with a subsequent extension onto $\mathbb Z_p$ (see \cite{R}). Then \cite{R} $S_2\varphi =\varphi \underline{*}1$ where 1 is the function equal identically to 1.

Let us describe invariant subspaces (we consider only closed subspaces) of the operator $S$. It will be convenient to use the model $S_1$.

\medskip
\begin{teo}
Every invariant subspace of the operator $S_1$ has the form $PH(\mathbb A_p)$, $P\in \mathbb C_p[x]$. The lattice $\operatorname{Lat} S_1$ of invariant subspaces is isomorphic to the division lattice of monic polynomials -- if $P',P''\in \mathbb C_p[x]$ are monic, then $P'H(\mathbb A_p)\subseteq P''H(\mathbb A_p)$ if and only if $P''$ divides $P'$. The commutant of $S_1$ coincides with the set of all operators of multiplication by functions from $H(\mathbb A_p)$, and also with the closure, in the uniform operator topology, of the set of all polynomials in the operator $S_1$.
\end{teo}

\medskip
{\it Proof}. It follows from the properties of $H(\mathbb A_p)$ as a Banach algebra that an invariant subspace of the operator $S_1$ is an ideal of $H(\mathbb A_p)$. Each ideal of $H(\mathbb A_p)$ has the form $PH(\mathbb A_p)$ where $P$ is a polynomial (see Corollary 5.5.11 in \cite{RC}). On the other hand, each of the above ideals is closed (thus forming an invariant subspace). Indeed, suppose that $P(z)=c(z-z_1)\cdots (z-z_k)$, $c\in \mathbb C_p$, $z_1,\ldots z_k\in \mathbb A_p$, and $Pf_n\to g$ uniformly. Here $f_n,g\in H(\mathbb A_p)$. The function $g$ admits the decomposition $g=P'h$ where $P'$ is a polynomial, while $h\in H(\mathbb A_p)$ has no zeroes (see the proof of the same corollary in \cite{RC} or Corollary 14.5 in \cite{E}). Then $z_1,\ldots z_k$ are among the zeroes of $P'$, so that $P$ divides $P'$ and $g\in PH(\mathbb A_p)$.

Thus, we have found a one-to-one correspondence between invariant subspaces of $S_1$ and ideals of the required form. The assertion regarding $\operatorname{Lat} S_1$ is obvious.

Suppose that $AS_1=S_1A$ where $A$ is a bounded linear operator on $H(\mathbb A_p)$. Consider an orthonormal basis in $H(\mathbb A_p)$, of the form $e_n(z)=z^n$, $n=0,1,2,\ldots$. Let $\varphi =Ae_0$. We have $\varphi \in H(\mathbb A_p)$,
$$
\varphi e_n=e_n\varphi =S_1^n\varphi =S_1^nAe_0=AS_1^ne_0=Ae_n.
$$
Therefore for any polynomial $\pi \in \mathbb C_p[z]$, we obtain the relation $A\pi =\varphi \pi$. Since polynomials are dense in $H(\mathbb A_p)$, we find that $Af=\varphi f$ for any $f\in H(\mathbb A_p)$. Approximating the function $\varphi$ uniformly by polynomials we come to the desired approximation property of the operator $A$. $\qquad \Box$

\bigskip
{\it Remarks}. 1) Classically, the Volterra integration operator on $L^2(0,1)$ has a linearly ordered lattice of invariant subspaces consisting of the subspaces of functions vanishing almost everywhere on $[0,\lambda ]$, $0<\lambda <1$. In our case, the operator $S_2$ has a linearly ordered sublattice of invariant subspaces
$$
X_N=\{ \varphi \in C(\mathbb Z_p,\mathbb C_p):\ \varphi (j)=0\text{ for $j\in \mathbb Z_+$, $0\le j\le N$}\}.
$$
However this sublattice does not exhaust the full lattice of invariant subspaces, as it can be seen from the above description given for $S_1$ -- the lattice $\operatorname{Lat} S_1$ is not linearly ordered.

2) Classically, the commutant of the unilateral shift on $H^2$ is known as the set of analytic Toeplitz operators (see Chapters 14 and 20 in \cite{H}) and coincides with the set of operators of multiplication by functions from $H^\infty$. The different situation in the non-Archimedean case is a consequence of the fact that $H(\mathbb A_p)$ is an algebra whereas $H^2$ is not.

\section{Backward Shift}
The backward shift operator $T:\ c_0\to c_0$ is defined as
\begin{equation}
T(x_0,x_1,x_2,\ldots )=(x_1,x_2,\ldots ),
\end{equation}
for each sequence $(x_n)\in c_0$. Sometimes it is convenient to use the same formulas (1) and (2) to define the extended shift operators acting on $l^\infty$. We denote them respectively by $\tilde{S}$ and $\tilde{T}$. Then it is easy to find the conjugate operators:
$$
S^*=\tilde{T},\quad T^*=\tilde{S}.
$$

Using the invariant subspaces of the operator $S$, we can construct some invariant subspaces of the operator $T$. Namely, let $N\subset c_0$ be invariant for $S$. Consider the subspace
$$
N^0=\{ x\in l^\infty :\ \langle x,y\rangle =0 \text{ for all $y\in N$}\}.
$$
For each $x\in N^0$, $\langle \tilde{T}x,y\rangle =\langle x,Sy\rangle =0$ for all $y\in N$. Therefore $\tilde{T}x\in N^0$. Thus, $N^0\cap c_0$ is an invariant subspace of the operator $T$.

Let us construct several functional models of the operator $T$.

\medskip
1) Let $T_1$ be an operator on $H(\mathbb A_p)$ of the form
$$
(T_1f)(z)=\frac{f(z)-f(0)}z,\quad z\in \mathbb A_p.
$$
If $f(z)=\sum\limits_{n=0}^\infty a_nz^n$, then $(T_1f)(z)=\sum\limits_{n=1}^\infty a_nz^{n-1}=\sum\limits_{n=0}^\infty a_{n+1}z^n$, that is under the isomorphism $H(\mathbb A_p)\to c_0$ the operator $T_1$ is transformed into $T$.

2) Let $T_2$ be an operator on $C(\mathbb Z_p,\mathbb C_p)$ of the form
$$
(T_2\varphi )(x)=\varphi (x+1)-\varphi (x),\quad x\in \mathbb Z_p.
$$
It was shown in \cite{K1} that $T_2$ can be interpreted as the annihilation operator in a $p$-adic representation of the canonical commutation relations of quantum mechanics. The operator $T_2$ acts on the Mahler polynomials as follows:
$$
T_2P_n=P_{n-1},\quad n\ge 1; \quad T_2P_0=0.
$$
Therefore under the isomorphism $\sum\limits_{n=0}^\infty c_nP_n\mapsto (c_0,c_1,c_2,\ldots )$ between the spaces $C(\mathbb Z_p,\mathbb C_p)$ and $c_0$, the operator $T_2$ is transformed into $T$.

Note \cite{K1} that the point spectrum of the operator $T_2$ coincides with the open ball $\{ \lambda \in \mathbb C_p:\ |\lambda |_p<1\}$; the eigenfunctions (``$p$-adic coherent states'') thave the form $\varphi_\lambda (x)=(1+\lambda )^x$, $x\in \mathbb Z_p$, so that $T_2\varphi_\lambda =\lambda \varphi_\lambda$.

3) Let us consider the Banach space $X_3$ of power series
\begin{equation}
g(z)=\sum\limits_{n=0}^\infty b_n\frac{z^n}{n!},\quad b_n\in \mathbb C_p,\quad |b_n|_p\to 0,
\end{equation}
with the norm $\|g\|_3=\sup\limits_{n\ge 0}|b_n|_p$. Since $|n!|_p=p^{-\nu (n)}$ where $\nu (n)\sim \frac{n}{p-1}$ (see Theorem 25.6 in \cite{S}), the space $X_3$ may be interpreted as the space of analytic elements on the closed ball $\left\{ z\in \mathbb C_p:\ |z|_p\le p^{-\frac1{p-1}}\right\}$. If we define an operator $T_3$ as the differentiation $(T_3g)(z)=g'(z)$, $g\in X_3$, then under the isomorphism $X_3\to c_0$, $\sum\limits_{n=0}^\infty b_n\dfrac{z^n}{n!}\mapsto (b_n)\in c_0$, the operator $T_3$ is transformed into $T$.

\medskip
The next result has a well-known classical counterpart (see Chapter 2 in \cite{N}).

\medskip
\begin{teo}
The operator $T$ has a dense set of cyclic vectors.
\end{teo}

\medskip
{\it Proof}. Let a vector $x=(x_0,x_1,x_2,\ldots )\in c_0$ be such that $x_k\ne 0$ for large values of $k$ (say $k\ge k_0$), and
\begin{equation}
\max\limits_{j\ge 1}\frac{|x_{k+j}|_p}{|x_k|_p}\longrightarrow 0,\quad \text{as $k\to \infty$}.
\end{equation}
Then $x$ is a cyclic vector.

Indeed, denote $e_n=(0,\ldots ,0,1,0,\ldots )$ where 1 is on the $n$-th place. Then
$$
\left\| \frac1{x_k}T^kx-e_0\right\| =\left\| \frac1{x_k}(x_k,x_{k+1},\ldots )
-e_0\right\| =\max\limits_{j\ge 1}\frac{|x_{k+j}|_p}{|x_k|_p}\longrightarrow 0,\quad k\to \infty,
$$
that is the vector $e_0$ belongs to the closed linear span $X'$ of the set $\{ x,Tx,T^2x,\ldots \}$.

For $k\ge k_0$, we have also
$$
T^{k-1}x-x_{k-1}e_0=(x_{k-1},x_k,\ldots )-(x_{k-1},0,\ldots )=(0,x_k,x_{k+1},\ldots )\in X',
$$
and
$$
\left\| \frac1{x_k}\left( T^{k-1}x-x_{k-1}e_0\right) -e_1\right\| =\max\limits_{j\ge 1}\frac{|x_{k+j}|_p}{|x_k|_p}\longrightarrow 0,\quad k\to \infty,
$$
so that $e_1\in X'$. Repeating this reasoning we show by induction that $e_n\in X'$ for each $n$, so that $X'=X$ and $x$ is a cyclic vector.

The condition (4) is fulfilled, for example, if $|x_k|_p=p^{-p^k}$ for $k\ge k_0$. Let $y\in c_0$. For any given $\varepsilon >0$, we can choose $k_0$ in such a way that the condition (4) holds, and $0<|x_k|_p<\varepsilon$, and also $|y_k|_p<\varepsilon$, as $k\ge k_0$. Set
$$
\tilde{y}=(y_0,y_1,\ldots ,y_{k_0-1},x_{k_0},x_{k_0+1},\ldots ).
$$
As we have proved, $\tilde{y}$ is a cyclic vector. On the other hand,
$$
\|y-\tilde{y}\|=\max\limits_{j\ge k_0}|y_j-x_j|_p<\varepsilon .
\qquad \Box
$$

\bigskip
More generally, if $E$ is a Banach space over $\mathbb C_p$, we can consider a backward shift operator $T_E$ on the space $c_0(E)$ of sequences $x=(x_0,x_1,x_2,\ldots )$, $x_j\in E$, such that $\|x_j\|_E\to 0$, as $j\to \infty$, with the norm $\|x\|=\sup\limits_{j\ge 0}\|x_j\|_E$. As before, $T_Ex=(x_1,x_2,\ldots )$. Obviously, $T_E$ is a contraction, that is $\|T_E\| \le 1$.

As in the classical case (see Problem 121 in \cite{H}), the operator $T_E$ has the following universality property.

\medskip
\begin{teo}
Let $A$ be a contraction on $E$, such that $\|A^nu\|_E\to 0$ for any $u\in E$. Then there exist a subspace $Y\subset c_0(E)$, invariant for the operator $T_E$, and a surjective isometry $W:\ E\to Y$, such that
\begin{equation}
A=W^{-1}T_E^YW,
\end{equation}
where $T_E^Y$ is the restriction of $T_E$ to $Y$.
\end{teo}

\medskip
{\it Proof}. Define $W:\ E\to c_0(E)$ setting
$$
Wu=(u,Au,A^2u,\ldots ),\quad u\in E.
$$
It is clear that $W$ is an isometry. Let $Y$ be the range of the operator $W$ in $c_0(E)$. Since $W$ is an isometry, $Y$ is closed. We have
$$
WAu=(Au,A^2u,\ldots )=T_E(u,Au,A^2u,\ldots )=T_EWu,
$$
so that $Au=W^{-1}T_EWu$. The equality $T_EWu=WAu$ also shows that $Y$ is invariant for $T_E$, and we come to the equality (5). $\qquad \Box$

\section*{ACKNOWLEDGEMENT}
This work was supported in part by the Ukrainian Foundation
for Fundamental Research, Grant 29.1/003.

\end{document}